\documentclass[a4paper,12pt]{article}

\usepackage[latin1]{inputenc}
\usepackage{amsmath}
\usepackage{amsfonts}
\usepackage{amssymb}
\usepackage{amsthm}
\usepackage{fontenc}
\usepackage[pdftex]{graphicx}
\usepackage[left=2.5cm,right=2.5cm,top=2cm,bottom=2cm]{geometry}
\usepackage{cite}
\usepackage{makeidx}
\usepackage{graphicx}
\usepackage{tikz}
\usepackage{pgf}

\newcommand{\R}{\mathbb{R}}

\newcommand{\Li}{{\mbox{Lip}}}

\newcommand{\diam}{\mathrm{diam}}
\newcommand{\dist}{\mathrm{dist}}

\newtheorem{definition}{Definition}
\newtheorem{theorem}{Theorem}
\newtheorem{lemma}{Lemma}

\newcommand{\E}{{\bf E}}

\begin{document}

\title{A higher Dimensional Marcinkiewicz Exponent and the Riemann Boundary Value Problems for Polymonogenic Functions on Fractals Domains}

\author{Carlos Daniel Tamayo Castro$^{a}$, Juan Bory Reyes$^{b}$}
\date { \small{$^{a}$Instituto de Matem\'aticas. Universidad Nacional Aut\'onoma de M\'exico, Mexico City, Mexico.\\
		cdtamayoc@comunidad.unam.mx\\
		$^{b}$SEPI-ESIME-Zacatenco, Instituto Polit\'ecnico Nacional, Mexico City, Mexico\\
	juanboryreyes@yahoo.com}
}

\maketitle

\begin{abstract} 
We use a high-dimensional version of the Marcinkiewicz exponent, a metric characteristic for non-rectifiable plane curves, to present a direct application to the solution of some kind of Riemann boundary value problems on fractal domains of Euclidean space $\mathbb{R}^{n+1}, n\geq2$ for Clifford algebra-valued polymonogenic functions with boundary data in classes of higher order Lipschitz functions. Sufficient conditions to guarantee the existence and uniqueness of solution to the problems are proved. To illustrate the delicate nature of this theory we described a class of hypersurfaces where the results are more refined than those that exist in literature.
\vspace{0.3cm}

\small{
  }  
\end{abstract}

\section{Introduction}
Riemann boundary value problems (RBVP for short) for analytic functions in a bounded Jordan domain of the complex plane are widely discussed and applied to many branches of mathematics, physics and engineering. Along classical lines, the Cauchy type integral is used as the main tool in the treatment of these boundary value problems, see \cite{Gajov, Lu, Mu} for more details.

Although the Cauchy type integral loses its meaning on fractal curves, the RBVP is equally valid. Pioneering works in the field were given by Kats in \cite{BKats83}. In particular, solvability conditions involving the H\"older exponent of the data function associated with the problem and the upper Minkowski dimension of the boundary are established. For a recent account of this approach we refer the reader to \cite{AbreuBlaya2011,Kats2009,Kats2014}.

The Marcinkiewicz exponent was introduced in \cite{BDKats16, 1DKats16, 2DKats16}, once again restricting the discussion to the complex analysis context. The use of this new metric characteristic of the boundary of the domains, considerably sharp the solvability conditions in the mentioned above study of the RBVPs.

Clifford analysis \cite{BDS, GHS08} offers an elegant generalization of analytic functions from the complex plane to higher-dimensions. Monogenic functions on Euclidean space are basic to Clifford analysis, they are defined as smooth functions with values in the corresponding Clifford algebra, which are null solutions to a generalized Cauchy-Riemann operator. In view of the fact that the generalized Cauchy-Riemann operator factorizes the higher Laplace operator, monogenic functions are also harmonic.

As was early pointed out in \cite[pp 22, 24]{Ryan96}, significant obstacles exist when trying to give a thorough treatment to the RBVP for monogenic functions. These are influenced by the fact that the product of two monogenic functions is not necessarily a monogenic. This is due to the non-commutativity properties of the Clifford algebra. This explains why an explicit solution to these boundary value problems has been found only for the jump problem and some slight modifications, see \cite{AB2001,Brs2001} and the references given there.

Almost all uniqueness and existence theorems for the solutions of RBVPs for analytic functions on the complex plane involving the Minkowski dimension of the boundary can be set in the context of Clifford analysis, see for instance \cite{AbreuBlaya2013, APB2007, TamayoCastro2022}. 

In \cite{AAB2015, Cerejeiras2012, YudeYuan2009, MinKaDa2012,LeDu2012} a kind of RBVPs have been studied for polymonogenic functions, i.e., null solutions to iterated generalized Cauchy-Riemann operator.

However, to the best of our knowledge, a full description of the relations between the high-dimensional version of the Marcinkiewicz exponent, introduced in \cite{TamayoCastro2022} and the RBVPs in fractal domains of Euclidean spaces for polymonogenic functions with boundary data in classes of higher order Lipschitz functions still remains open.

This paper aims to obtain solvability conditions to RBVPs for Clifford algebra-valued polymonogenic functions involving the high-dimensional Marcinkiewicz exponent of the hypersurface. Besides, to show that this condition strictly improves those requiring the Minkowski dimension, we need to describe an uncountable class of hypersurfaces, such that the relationship between its Minkowski dimension and Marcinkiewicz exponent is with strict inequality.

The content of this paper is structured as follows: Section \ref{Preliminars} presents the basic notions and terminology of the theory of monogenic functions. Moreover, we give a brief review of the definitions of higher order Lipschitz functions, fractal dimensions and Marcinkiewicz exponent. Also, a Whitney type extension theorem is stated. In Section \ref{BoundaryValueProblem}, we obtain solvability and uniqueness conditions to a kind of RBVPs for Clifford algebra-valued polymonogenic functions with data higher order Lipschitz functions. Section \ref{Example} provides a description of a class of hypersurfaces, where the conditions obtained in Section \ref{BoundaryValueProblem} improve those involving the Minkowski dimension that exist in literature.

\section{Preliminaries and Notations}\label{Preliminars}
This section presents the essential background needed for developing the results in the subsequent sections. It is divided into three subsections, each devoted to a fundamental component of the strategy for proving the results. The first subsection contains some basic facts
about Clifford analysis. In second subsection, we give a brief exposition of the data classes of functions required in the problems considered, and a corresponding Whitney-type extension theorem is presented. Finally, we introduce the fractal dimensions and characteristics of the boundaries that will be treated in the work.      
\subsection{Clifford Algebras and Monogenic Functions}
In this section, we recall some basic facts about Clifford analysis which will be needed in the sequel. It could be seen in more detail in the literature, for instance \cite{BDS, GHS08}.
\begin{definition}\label{CliffAlg}
	We consider the Clifford algebra associated with $\mathbb{R}^{n}$ and endowed with the usual Euclidean metric, as the minimal extension of $\mathbb{R}^{n}$ to a unitary, associative algebra $\mathcal{C}\ell(n)$ over the real numbers, satisfying 
	\begin{eqnarray*}
		x^{2} = -\arrowvert x\arrowvert^{2},   
	\end{eqnarray*}
	for any  $x \in \mathbb{R}^{n}$.
\end{definition}
It thus follows that, the standard orthonormal basis of $\mathbb{R}^{n}$, denoted by $\{e_{j}\}_{j = 1}^{n}$, is subject to the
basic multiplication rules
\begin{displaymath}
	e_{i}e_{j} + e_{j}e_{i} = -2\delta_{ij},
\end{displaymath}
were $\delta_{ij}$ is the Kronecker delta.

For arbitrary $a \in \mathcal{C}\ell(n)$ we have $a = \sum\limits_{A \subseteq N}a_{A}e_{A}$, $N = \{1, \ldots, n \}$, $a_{A} \in \mathbb{R}$ where $e_{\emptyset} = e_{0} = 1$, $e_{\{ j\}} = e_{j}$ and $e_{A} = e_{\beta_{1}}\cdots e_{\beta_{k}}$ for $A = \{\beta_{1}, \ldots, \beta_{k} \}$ with  $\beta_{j} \in \{1, \ldots, n \}$ and $\beta_{1} < \ldots < \beta_{k}$.

An important subspace of $\mathcal{C}\ell(n)$ is the so-called space of paravectors and its elements have the form $x = \sum_{j = 0}^{n}x^{j}e_{j},$. Each $x = (x^{0}, x^{1}, \ldots, x^{n}) \in \mathbb{R}^{n + 1}$ will be identified with a paravector.

The algebra norm of arbitrary $a$ is defined by $\arrowvert a\arrowvert := \left(\sum_{A}a_{A}^{2}\right)^{\frac{1}{2}}$ and $\mathcal{C}\ell(n)$ becomes an Euclidean space.

The conjugation in $\mathcal{C}\ell(n)$ is defined as the anti-involution $a\mapsto\overline{a}:= \sum_{A}a_{A}\overline{e}_{A}$ such that
\begin{displaymath}
	\overline{e}_{A} := (-1)^{k}e_{\beta_{k}} \cdots e_{\beta_{2}}e_{\beta_{1}}.
\end{displaymath}
For every paravector $x$ the relation
\begin{displaymath}
	x\overline{x} = \overline{x}x = |x|^{2}
\end{displaymath}
holds.

An $\mathcal{C}\ell(n)-$valued function $u$ defined over $\Omega\subset\mathbb{R}^{n + 1}$ has the representation
\begin{displaymath}
	u(x) = \sum_{A}u_{A}(x)e_{A},
\end{displaymath}
where $u_{A}$ are $\mathbb{R}$-valued components.

Unless otherwise states we assume that $u$ to be $\mathcal{C}\ell(n)-$valued. Properties such as continuity or differentiability have to be understood component-wise.

Clifford analysis is mainly centered around the concept of monogenic functions which are null solutions of the generalized Cauchy-Riemann operator in $\mathbb{R}^{n + 1}$ defined by
\begin{displaymath}
	\mathcal{D} := \sum^{n}_{j = 0}e_{j}\frac{\partial}{\partial x_{j}}.
\end{displaymath}
Let $\Omega\subset\mathbb{R}^{n + 1}$ be an open set and $u \in C^{1}(\Omega)$ then $u$ will be called left (respectively right) monogenic in $\Omega$ if $\mathcal{D}u = 0$ (respectively $u\mathcal{D} = 0$) in $\Omega$.

The fundamental solution for the Cauchy-Riemann operator is  
\begin{displaymath}\label{Ch1EqFundSolCauchRiem}
	\begin{array}{cc}
		E(x) = \dfrac{1}{\sigma_{n + 1}}\dfrac{\overline{x}}{|x|^{n + 1}}, &  x\in\mathbb{R}^{n + 1}\setminus \{0\}
	\end{array} 
\end{displaymath}
where $\sigma_{n + 1}$ stands for the hypersurface area of the unit sphere in $\mathbb{R}^{n + 1}$. 

Polymonogenic functions $f \in C^{k}(\Omega)$ of order $k$, emerge as the solutions of the iterated action of the operator $\mathcal{D}$
\begin{displaymath}
	\mathcal{D}^{k}f = 0, 
\end{displaymath}
in $\Omega \subset \mathbb{R}^{n + 1}$. This notion goes back as far as \cite{Brgehr1999,Brackx1976,DelangheBrackx1978}.  

For arbitrary $k \in N$, let us introduce 
\begin{displaymath}
	E^{k}(x) =	\dfrac{1}{\sigma_{n+1}}\dfrac{\overline{x}(\overline{x} + x)^{k-1}}{2^{k-1}(k-1)!\lvert x\rvert^{n+1}}.
\end{displaymath}
Note that for $k = 1$ 
\begin{eqnarray*}
	E^{1}(x) = E(x), 	
\end{eqnarray*}
as it is easy to check.

A direct computation shows that
\begin{eqnarray*}
	\mathcal{D}E^{k}(x) = E^{k-1}(x).
\end{eqnarray*}
Consequently, by decreasing induction on $k$, 
\begin{eqnarray*}
	\mathcal{D}^{k}E^{k}(x) = \mathcal{D}E^{1}(x) = 0, &  x\in\mathbb{R}^{n + 1}\setminus \{0\}.
\end{eqnarray*} 
Let us now introduce the Teodorescu transform related to the theory of monogenic functions, see \cite{GHS08} for more details.
\begin{definition}\label{Ch1DefTeodTransf}
	Let $\Omega\subset \mathbb{R}^{n+1}$ be a domain and let $u \in C^{1}(\overline{\Omega})$, the so-called Teodorescu transform is defined by
	\begin{eqnarray*}
		(T_{\Omega}u)(x) := -\int\limits_{\Omega}E(y-x)u(y)dV(y), & x \in \mathbb{R}^{n + 1}, 
	\end{eqnarray*}
	where $dV(y)$ is the volume element. 
\end{definition}  
Sufficient conditions for the H\"older continuity of $T_{\Omega}u$ are established by our next theorem.
\begin{theorem}\label{PropTeo1}
	Let $\Omega\subset \mathbb{R}^{n+1}$ be a domain and let $u \in L^{p}(\Omega)$ for $p > n + 1$. Then
	\begin{itemize}
		\item The integral $(T_{\Omega}u)(x)$ exists in the whole $\mathbb{R}^{n + 1}$ and tends to zero as $\arrowvert x\arrowvert \rightarrow \infty$. Besides, $T_{\Omega}u$ is a monogenic function in $\mathbb{R}^{n + 1}\setminus \overline{\Omega}$. 
		\item For $x, y \in \mathbb{R}^{n + 1}$, and $x \neq y$, we have 
		\begin{displaymath}
			\arrowvert(T_{\Omega}u)(x) - (T_{G}u)(y)\arrowvert \leq C_{2}(\Omega, p, n)\|u\|_{p}\arrowvert x - y\arrowvert^{\frac{p - n - 1}{p}}.
		\end{displaymath}
	\end{itemize}
\end{theorem}
Conditions for the derivability of $T_{\Omega}u$ on $\Omega$ are given as follows
\begin{theorem}\label{PropTeo2}
	Let $\Omega\subset \mathbb{R}^{n+1}$ be as stated above. Then $T_{\Omega}u$ is a differentiable function for every $x \in \Omega$ with 
	\begin{displaymath}
		\frac{\partial}{\partial x_{j}}(T_{\Omega}u)(x) = -\int\limits_{\Omega}\frac{\partial}{\partial x_{j}}[E_{n}(y - x)]u(y)dV(y) + \overline{e_{j}}\dfrac{u(x)}{n + 1}.
	\end{displaymath}
	In particular, the identity 
	\begin{eqnarray*}
		\mathcal{D}(T_{\Omega}u)(x) = u(x),\ \ x \in \Omega.
	\end{eqnarray*}
	holds.
\end{theorem}
Now, we introduce a kind of polymonogenic Teodorescu transform following \cite{Brgehr1999}. Due to its good properties, this integral operator will play an essential role in the strategies described below. For a recent account of the subject we refer the reader to \cite{Liping2016}.
\begin{definition}
	Let $\Omega \subset \mathbb{R}^{n+1}$ be a bounded domain and $u \in L^{1}(\overline{\Omega})$. Then for $k \in \mathbb{N}$ we define the $k$-polymonogenic Teodorescu transform as follow
	\begin{displaymath}
		T^{k}_{\Omega}u(x) := (-1)^{k}\int\limits_{\Omega}E^{k}(y - x)u(y)dV(y),
	\end{displaymath}
	where $dV(y)$ is the volume element.
\end{definition}
By derivability properties of $T^{k}_{\Omega}u(x)$, the following equality is achieved
\begin{eqnarray*}
	\mathcal{D}T^{k}_{\Omega}u = T^{k -1}_{\Omega}u, & k \geq 2.
\end{eqnarray*}
Therefore, decreasing induction on $k$ combining with Theorems \ref{PropTeo1} and \ref{PropTeo2} gives
\begin{equation}\label{Ch4EqDkTk}
	\mathcal{D}^{k}T^{k}_{\Omega}u = \mathcal{D}T^{1}_{\Omega}u = \left\lbrace \begin{array}{ccc}
		u,	& in & \Omega \\
		0,	& in & \mathbb{R}^{n+1}\setminus \overline{\Omega}.
	\end{array} \right. 
\end{equation}

\subsection{Function Classes and a kind of a Whitney Type Theorem}
In order to state the main problems that we address, it is necessary to define the suitable classes where the data functions will be considered. These are the higher order Lipschitz classes, on the basis of which some Whitney-type extension theorems will be presented.

As a first step, we shall recall the class of $p$-integrables functions, see \cite{GHS08}.  
\begin{definition}
	Let  $\Omega \subset \mathbb{R}^{n + 1}$ be a domain and let $u: \Omega \mapsto \mathcal{C}\ell(n)$. $L^{p}(\Omega), 0 < p < \infty$ denotes the space of all equivalence classes of Lebesgue measurable functions equal almost everywhere, such that 
	\begin{displaymath}
		\|u\|_{p} < \infty.
	\end{displaymath}
	where
	\begin{displaymath}
		\|u\|_{p} := \left(\int\limits_{\Omega}|u|^{p}dV \right)^{\frac{1}{p}},
	\end{displaymath}	
\end{definition}
Let $\E$ be a closed subset of $\R^{n + 1}$. We write  $j=(j_0, j_1, \cdots, j_n)$ a n-dimensional multi-index of order $|j| = j_0 + j_1 +\cdots+j_n$, where $j_0,  j_1, \dots,j_n$ are non-negative integers. In addition, we have, $j! = j_0!j_1! \cdots j_n!$ and $x^{j} = x_{j_{0}}x_{1}^{j_{1}}\cdots x_{n}^{j_{n}}$. 
\begin{definition}
	Let $0 < \nu \leq 1$. We shall say that a real-valued function $f$, defined in $\E$, belongs to $\Li(\E,k + \nu)$ if there exist real-valued bounded functions  $f^{(j)}$, $0<|j|\le k$, defined on $\E$, with $f^{(0)}=f$, and so that
	\begin{equation}\label{L1}
		f^{(j)}(x)=\sum_{|j+l|\le k}\frac{f^{(j+l)}(y)}{l!}(x-y)^l+R_j(x,y),\,\,x,y\in\E,
	\end{equation}
	where
	\begin{equation}\label{L2}
		|f^{(j)}(x)|\le M,\,\,\,|R_j(x,y)|\le M|x-y|^{k + \nu -|j|},\,\,x,y\in\E, |j|\le k,
	\end{equation}
	being $M$ a positive constant.
\end{definition}
The space of all these functions is named higher order Lipschitz spaces, see for instance \cite{St}. When $k = 1$ and $\E$ is a compact set,  these functions reduce to the classical H\"older continuous functions.

A norm in $\Li(\E,k + \nu)$ is defined as the smallest $M$ satisfying (\ref{L2}). In \cite{TamayoCastro2018} was proved that $\Li(\E,k + \nu)$ endowed with this norm is a Banach space. Also, conditions for continuous and compact embeddings of generalized higher-order Lipschitz classes on a compact subset of Euclidean spaces were obtained, showing that these spaces are not only a generalization but also a refinement of the classical Lipschitz classes. 

In general, an element of $\Li(\E,k + \nu)$ should be interpreted as a collection 
\begin{displaymath}
	\left\lbrace f^{(j)}:\E\mapsto\R,\,|j|\le k\right\rbrace.	
\end{displaymath}
In order to present a kind of a Whitney type extension theorem for polymonogenic functions, let us recall the following classical theorem by Whitney, see \cite[pag. 177]{St}. 
\begin{theorem}\label{Ch1WhitneyExtLip}
	Let $\textbf{E} \subset \mathbb{R}^{n+1}$ be a closed set and let $f \in \Li(\textbf{E}, k + \nu)$ with values in $\mathbb{R}$. Then, there exists a $\mathbb{R}$-valued function $\widetilde{f} \in \Li(\mathbb{R}^{n+1}, k + \nu)$ satisfying 
	\begin{itemize}
		\item $\partial^{(j)}\widetilde{f}\arrowvert_{E} = f^{(j)}$,  
		\item $\widetilde{f} \in C^{\infty}(\mathbb{R}^{n+1}\setminus \textbf{E})$,
		\item $\arrowvert \partial^{(j)}\widetilde{f}(x) \arrowvert \leq C\dist(x, \textbf{E})^{\nu - 1}$ for $\arrowvert j\arrowvert = k+1$ and $x \in \mathbb{R}^{n+1}\setminus \textbf{E}$.
	\end{itemize}
\end{theorem}
Here and subsequently the symbol 
\begin{displaymath}
	\partial^{(j)}:=\frac{\partial^{|j|}}{\partial x_0^{j_0}\partial x_1^{j_1} \cdots \partial x_k^{j_k}},	
\end{displaymath}
stands for the higher-order partial derivatives.

Let $f \in \Li(\mathcal{S},k-1 + \nu)$ be a $\mathcal{C}\ell(n)$-valued function, interpreted as the collection 
$$\left\lbrace f^{(j)}:\mathcal{S}\mapsto\mathcal{C}\ell(n),\,|j|\le k-1\right\rbrace $$ 
with $f^{(0)} = f$ satisfying (\ref{L1}) and (\ref{L2}). 

In order to present a suitable version of Whitney extension theorem for $\mathcal{C}\ell(n)$-valued function, in \cite{AAB2015} are constructed the following functions,
\begin{eqnarray}\label{Ch1EqFunctFiPoly}
	\mathbf{f}^{(i)} = \sum_{r_{1},\cdots,r_{i} = 0}^{n}e_{r_{1}}\cdots e_{r_{i}}f^{\textbf{1}_{r_{1}}+\cdots+\textbf{1}_{r_{i}}}, & i = 0, 1, \dots, k - 1.
\end{eqnarray}
Here $\textbf{1}_{r_{i}}$ is the multi-index $(j_{0}, j_{1},\cdots, j_{n})$ with 
\begin{equation}\label{Ch1EqEpecialMultiFi}
	j_{p} = \left\lbrace \begin{array}{cc}
		1, & p = r_{i}	\\
		0, & p \neq r_{i}.	
	\end{array}\right. 
\end{equation}
We should note that the functions $\mathbf{f}^{(i)}$ are an appropriate arranged of every function $f^{(j)}$ with $\left| j \right| = i$. In addition, $\mathbf{f}^{(0)}$ = $f^{(0)}$ = $f$. 

Let us mention an important consequence of Theorem \ref{Ch1WhitneyExtLip}. This can be found in \cite{AAB2015}.
\begin{theorem}\label{Ch4WhitneyExtPoly}
	Let $\textbf{E} \subset \mathbb{R}^{n+1}$ be a closed set and let $f \in \Li(\textbf{E}, k-1 + \nu)$ with values in $\mathcal{C}\ell(n)$. Then, there exists a $\mathcal{C}\ell(n)$-valued function $\widetilde{f} \in \Li(\mathbb{R}^{n+1}, k-1 + \nu)$ satisfying 
	\begin{itemize}
		\item $\mathcal{D}^{i}\widetilde{f}\arrowvert_{E} = \mathbf{f}^{(i)}$, \, $i = 0, 1, \cdots k-1 $ 
		\item $\widetilde{f} \in C^{\infty}(\mathbb{R}^{n+1}\setminus \textbf{E})$,
		\item $\arrowvert \mathcal{D}^{k}\widetilde{f}(x) \arrowvert \leq C\dist(x, \textbf{E})^{\nu - 1}$ for $x \in \mathbb{R}^{n+1}\setminus \textbf{E}$.
	\end{itemize}
\end{theorem}

\subsection{Fractal Dimensions and Characteristics}\label{SubSecFracDim}
To make the presentation self-contained, we review some basic ideas about fractal dimensions of sets that will be required to work with domains with fractal boundaries. For a deeper discussion of this topic, we refer the reader to \cite{Fal, Mandelbrot, Matt}.

Let $\textbf{E} \in \mathbb{R}^{n+1}$ a non-empty set. For any $s \geq 0$ and $\delta > 0$, $\mathcal{H}_{\delta}^{s}(\textbf{E})$  is defined as, 
\begin{displaymath}
	\mathcal{H}_{\delta}^{s}(\textbf{E}) := \inf\{\sum_{i = 1}^{\infty}\diam(U_{i})^{s}: \{U_{i}\} \ is \ a  \ \delta-covering \ of \ \textbf{E}\},
\end{displaymath} 
where $\diam(U)$ is the diameter of the set $U$. The infimum here is taken over all countable $\delta$-coverings ${U_{i}}$ of $\textbf{E}$ for open or closed balls. Now, we can present the Hausdorff measure.
\begin{definition}
	The $s$-dimensional Hausdorff measure is defined by	the limit
	\begin{displaymath}
		\mathcal{H}^{s}(\textbf{E}) := \lim _{\delta\rightarrow 0}\mathcal{H}_{\delta}^{s}(\textbf{E}).
	\end{displaymath} 
\end{definition}
It can be shown that the $s$-dimensional Hausdorff measure of a set $\textbf{E}$ is almost always 0 or $\infty$. There is only one value of $s$ where the measure change between these two values. Therefore, it looks natural to define the Hausdorff dimension as this value. 
\begin{definition}\label{DefHausDim}
	The Hausdorff dimension of \textbf{E} is defined as
	\begin{displaymath}
		\dim _{H}\textbf{E} := \inf\{s \geq 0: \mathcal{H}^{s}(\textbf{E}) = 0\} = \sup\{s \geq 0: \mathcal{H}^{s}(\textbf{E}) = \infty\}.
	\end{displaymath}   
\end{definition}
The following theorem can be found in \cite{APB2007}.
\begin{theorem}\label{CorrDolz}
	Let $\Omega$ be a domain in $\mathbb{R}^{n+1}$ and $\textbf{E} \subset \Omega$ be a compact set. Let be $\mathcal{H}^{n + \mu}(\textbf{E}) = 0$ where $0 < \mu \leq 1$. If $u \in \Li(\Omega, \mu)$, and it is monogenic in $\Omega\setminus\textbf{E}$, then $u$ is also monogenic in $\Omega$.   
\end{theorem}  
The following definition of fractal set is due to Mandelbrot \cite{Mandelbrot}
\begin{definition}\label{DefFractMan}
	If an arbitrary set $\textbf{E} \subset \mathbb{R}^{n+1}$ with topological dimension $n$ has $\dim_{H}\textbf{E} > n$,  then $\textbf{E}$ is called a fractal set.
\end{definition}
The results presented in this paper are intended to deal with sets such that $\mathcal{H}^{n}(\textbf{E}) = \infty$ as well as fractals from Definition \ref{DefFractMan}. Throughout this paper, the expression \textquoteleft fractal domain' will always refer to a domain with a fractal boundary.

The Minkowski dimension is widely used when working with fractals. That is due to the fact that computations are easier than with other fractal dimensions.  We will only consider the upper Minkowski dimension. 
\begin{definition}\label{DefMinkDim}
	Let $\textbf{E}$ be a non-empty bounded subset of $\mathbb{R}^{n+1}$ and let $N_{\delta}(\textbf{E})$ be the smallest number of sets of diameter at most $\delta$, covering $\textbf{E}$. The upper Minkowski dimension of $\textbf{E}$ is defined as 
	\begin{displaymath}
		\overline{\dim}_{M}\textbf{E} := \limsup_{\delta \rightarrow 0}\dfrac{\log N_{\delta}(\textbf{E}) }{-\log \delta}. 
	\end{displaymath}
\end{definition}
We denote $\mathcal{M}_{0}$ the grid covering $\mathbb{R}^{n+1}$ which consists of $(n+1)$-dimensional cubes with vertices with integer coordinates and edges of length one. From $\mathcal{M}_{0}$ is generated the grid $\mathcal{M}_{k}$ by dividing every cube in $\mathcal{M}_{0}$ into $2^{(n + 1)k}$ different cubes with edges lengths $2^{-k}$. Let $N_{k}(\textbf{E})$ be the amount of cubes of the grid $\mathcal{M}_{k}$ which intersect $\textbf{E}$. Then as can be found in \cite{Fal} we have
\begin{equation}\label{Ch1EqEquiDefMinkDim}
	\overline{\dim}_{M}\textbf{E} = \limsup_{k \rightarrow \infty}\dfrac{\log N_{k}(\textbf{E}) }{k\log (2)}. 
\end{equation}
In \cite[pp 77]{Matt} is given the next theorem relating the Hausdorff and Minkowski dimensions.
\begin{theorem}
	For the bounded set $\textbf{E} \subset \mathbb{R}^{n+1}$ with topological dimension $n$, we have
	\begin{displaymath}
		n \leq \dim_{H}\textbf{E} \leq \overline{\dim}_{M}\textbf{E} \leq n + 1.
	\end{displaymath}	
\end{theorem}
The following new metric characteristics of a fractal set in $\mathbb{R}^{n + 1}$ are mainly included to keep the exposition as self-contained as possible. These can be found in \cite{TamayoCastro2022}.

Let $\mathcal{S}$ be a topologically compact hypersurface in $\mathbb{R}^{n +1}$, which bounds a Jordan domain $\Omega^{+}$. We write $\Omega^{-}$ for the unbounded complement. It is assumed $\mathcal{S}$ to be fractal.

Let $D \subset \mathbb{R}^{n + 1}$ be a bounded set, which does not touch the hypersurface $\mathcal{S}$. We will consider the integral
\begin{displaymath}
	I_{p}(D) = \int\limits_{D}\dfrac{dV(x)}{\dist^{p}(x, \mathcal{S})}.
\end{displaymath}
For completeness, we recall:
\begin{definition}\label{DefMarcExp}
	Let $\mathcal{S}$ be a topologically compact hypersurface which is the boundary of a Jordan domain in $\mathbb{R}^{n +1}$. We define the inner and outer Marcinkiewicz exponent of $\mathcal{S}$, respectively, as
	\begin{displaymath}
		\begin{array}{cc}
			\mathfrak{m}^{+}(\mathcal{S}) = \sup \{p: I_{p}(\Omega^{+}) < \infty\},  & \mathfrak{m}^{-}(\mathcal{S}) = \sup\{p: I_{p}(\Omega^{*}) < \infty\}, 
		\end{array} 
	\end{displaymath}
	and the (absolute) Marcinkiewicz exponent of $\mathcal{S}$ as,
	\begin{displaymath}
		\mathfrak{m}(\mathcal{S}) = \max\{\mathfrak{m}^{+}(\mathcal{S}), \mathfrak{m}^{-}(\mathcal{S})\}.
	\end{displaymath} 
\end{definition}
Here, the domain $\Omega^{*} := \Omega^{-}\bigcap\{x: \arrowvert x\arrowvert < r\}$, where $r$ is chosen in a way that $\mathcal{S}$ is completely contained inside the ball of radius $r$. We should note that the value of $\mathfrak{m}^{-}(\mathcal{S})$ does not depend on the selection of the radius $r$ when constructing $\Omega^{*}$, due to the fact that only the points closest to $\mathcal{S}$ influence the convergence of the integral $I_{p}(D)$.

The following theorem expresses the relationship between the Minkowski dimension with the Marcinkiewicz exponent, it was proved in \cite{2DKats16} and in a different way in \cite{TamayoCastro2022}. 
\begin{theorem}\label{TheoIneqRm}
	Let $\mathcal{S}$ be a topologically compact hypersurface which is the boundary of a Jordan domain in $\mathbb{R}^{n +1}$, then $\mathfrak{m}(\mathcal{S}) \geq n + 1 - \overline{\dim}_{M}(\mathcal{S})$.
\end{theorem}

\section{Riemann Boundary Value Problems}\label{BoundaryValueProblem}
In this section our main results concerning RBVPs for polymonogenic functions on fractal domains using the absolute Marcinkiewicz exponent are stated and proved. 

Let $f \in \Li(\mathcal{S}, k-1 + \nu)$ be a $\mathcal{C}\ell(n)-$valued function. We are first interested in the following boundary value problem: To find a polymonogenic function $\Phi$ of order k on $\R^{n+1}\setminus\mathcal{S}$ continuously extendable from $\Omega^{\pm}$ to $\mathcal{S}$ such that its boundary values $\Phi^{\pm}$ on $\mathcal{S}$ satisfy the following conditions
\begin{equation}\label{Ch4JumpProbPoly}
	\begin{array}{ccc}
		(\mathcal{D}^{i}\Phi(x))^{+} - (\mathcal{D}^{i}\Phi(x))^{-} = \mathbf{f}^{(i)} & x \in \mathcal{S} & 0 \leq i \leq k-1\\
		(\mathcal{D}^{i}\Phi(\infty))^{-} = 0 &  & 0 \leq i \leq k-1,
	\end{array}
\end{equation}
where the functions $\textbf{f}^{(i)}$ were defined in \eqref{Ch1EqFunctFiPoly}.

As a special case when  $k = 1$ the polymonogenic functions derive to monogenic functions. At the same time, the higher order Lipschitz class $\Li(\mathcal{S}, k-1 + \nu)$ becomes the standard Lipschitz class with exponent $\nu$. Consequently, problem \eqref{Ch4JumpProbPoly} reduces the classical jump problem for monogenic functions: 
\begin{equation}\label{JumpProbl}
	\begin{array}{ccc}
		\Phi^{+}(x) - \Phi^{-}(x) = f, & x \in \mathcal{S}, \\
		\Phi^{-}(\infty) = 0. &  
	\end{array}
\end{equation}
Hence, problem (\ref{Ch4JumpProbPoly}) generalizes problem (\ref{JumpProbl}) presented and studied in \cite{TamayoCastro2022}.

The following lemma was proved in \cite{AAB2015} and takes part of the proof of the upcoming theorem.
\begin{lemma}\label{Ch4LemLpPoly}
	Let $\Omega$ be a bounded domain of $\mathbb{R}^{n+1}$ and let $g \in L^{p}(\Omega)$ with $p > n+1$. Then,
	\begin{eqnarray*}
		\mathcal{D}^{i}T^{k}_{\Omega}g \in \Li(\mathbb{R}^{n+1}, \alpha), & i = 0, 1, \cdots, k - 1;
	\end{eqnarray*}
	with $0 < \alpha \leq \frac{p-n-1}{p}$.
\end{lemma}
The following theorem provides a sufficient solvability condition to problem (\ref{Ch4JumpProbPoly}) and generalizes and strengthens \cite[Theorem 9]{TamayoCastro2022}.
\begin{theorem}\label{Ch4TheoCondSolvPoly}
	If $f \in \Li(\mathcal{S}, k-1 + \nu)$, with
	\begin{eqnarray}\label{Ch4EqCondSolvPoly}
		\nu > 1 - \dfrac{\mathfrak{m}(\mathcal{S})}{n + 1},
	\end{eqnarray}
	and $k < n + 1$, then the problem (\ref{Ch4JumpProbPoly}) is solvable. 
\end{theorem}
\begin{proof}
	We need to show that the solution is given by  
	\begin{equation}\label{Ch4EqSolPoly+}
		\begin{array}{cc}
			\Phi(x) = \widetilde{f}(x)\chi^{+}(x) - (T^{k}_{\Omega^{+}}\mathcal{D}^{k} \widetilde{f})(x), & x \in \mathbb{R}^{n+1},
		\end{array}
	\end{equation}	
	when $\mathfrak{m}(\mathcal{S}) = \mathfrak{m}^{+}(\mathcal{S})$, or by
	\begin{eqnarray}\label{Ch4EqSolPoly*}
		\Phi(x) = -f^{*}(x)\chi^{*}(x) + (T^{k}_{\Omega^{*}}\mathcal{D}^{k} f^{*})(x), & x \in \mathbb{R}^{n+1},
	\end{eqnarray}
	when $\mathfrak{m}(\mathcal{S}) = \mathfrak{m}^{-}(\mathcal{S})$.
	
	Here, $\chi^{+}(x)$ and $\chi^{*}(x)$ are the characteristic functions of $\Omega^{+}$ and $\Omega^{*}$ respectively. Besides, $\widetilde{f}$ is the Whitney extension of $f$ and $f^{*} = \widetilde{f}\rho$ where $\rho$ is defined as follows. We will fix $r_{1}$ large enough such that $\mathcal{S}$ is entirely contained inside the ball $B_{1} = \{x: \arrowvert x\arrowvert < r_{1}\}$. We choose $r > r_{1}$, and define $B = \{z: \arrowvert x\arrowvert < r\}$. Due to the fact that the value of $\mathfrak{m}^{-}(\mathcal{S})$ do not depend on the selection of $r$ then, let be $\Omega^{*} = \Omega^{-}\bigcap B$. Thus, let  $\rho(x)$ be a real valued function in $C^{\infty}(\mathbb{R}^{n + 1})$ such that $0 \leq \rho(x) \leq 1$, equal to 0 outside of $B$, and equal to 1 over $B_{1}$. 	
	
	The same proof works for $\mathfrak{m}^{+}(\mathcal{S})$ and $\mathfrak{m}^{-}(\mathcal{S})$, we will consider the first case. 
	
	We must have that $\mathcal{D}^{k}\widetilde{f} \in$ L$^{p}(\Omega^{+})$ with $p > n + 1$, being $\widetilde{f}$ the Whitney extension of $f$. Theorem \ref{Ch4WhitneyExtPoly} now show that 
	\begin{equation*}
		\int\limits_{\Omega^{+}}\arrowvert\mathcal{D}^{k}\widetilde{f}(x)\arrowvert^{p}dV(x) \leq C\int\limits_{\Omega^{+}}\dfrac{dV(x)}{\dist(x, \mathcal{S})^{p(1 - \nu)}}.
	\end{equation*}
	The right-hand integral above converges for $p < \frac{\mathfrak{m}^{+}(\mathcal{S})}{1 - \nu}$, which is a direct consequence of Definition \ref{DefMarcExp},. Then the main requirement is  
	\begin{equation*}
		\nu > 1 - \dfrac{\mathfrak{m}^{+}(\mathcal{S})}{n + 1}.
	\end{equation*} 
	From (\ref{Ch4EqDkTk}) it follows that $\Phi$ is a polymonogenic function of order $k$ on $\R^{n+1}\setminus\mathcal{S}$. Combining Lemma \ref{Ch4LemLpPoly} with the fact that $\widetilde{f} \in \Li(\mathbb{R}^{n + 1}, k - 1 + \nu)$  we obtain that the functions $\mathcal{D}^{i}\Phi$, $i = 0, 1, \dots, k-1$; are continuous functions on $\overline{\Omega^{+}}$ and $\overline{\Omega^{-}}$. 
	
	Combining Lemma \ref{Ch4LemLpPoly} with Theorem \ref{Ch4WhitneyExtPoly}, we can conclude that the function $\Phi(x)$ satisfies the boundary condition over $\mathcal{S}$. Finally, as was stated in \cite{AAB2015} when $k < n+1$, we have that  $\mathcal{D}^{i}\Phi^{-}$ vanishes at infinity for every $i = 0, 1,\cdots, k-1$. A trivial verification shows that $\mathcal{D}^{k} f^{*} \in L^{p}(\Omega^{*})$ providing that $\nu > 1 - \frac{\mathfrak{m}^{-}(\mathcal{S})}{n + 1}$. This completes the proof.  
\end{proof}
We can also prove a sufficient condition for unicity. The next theorem is a generalization of \cite[Theorem 10]{TamayoCastro2022}.
\begin{theorem}\label{Ch4TheoUnicityCondPoly}
	Let be $f \in \Li(\mathcal{S}, k-1 + \nu)$ with $\nu > 1 - \dfrac{\mathfrak{m}(\mathcal{S})}{n + 1}$ and $k < n + 1$, let 
	\begin{equation*}
		\dim_{H}\mathcal{S} - n < \mu < 1 - \dfrac{(n+1)(1-\nu)}{\mathfrak{m}(\mathcal{S})}.
	\end{equation*}
	Then there is a unique solution $\Phi$ of the problem (\ref{Ch4JumpProbPoly}), such that $\mathcal{D}^{i}\Phi$  belongs to the classes $\Li(\overline{\Omega^{+}}, \mu)$ and $\Li(\overline{\Omega^{-}},\mu)$, for $i = 0,1,\cdots,k-1$.
\end{theorem}
\begin{proof}
	From Lemma \ref{Ch4LemLpPoly} and the proof of Theorem \ref{Ch4TheoCondSolvPoly} we deduce that the solution $\Phi$ to the problem (\ref{Ch4JumpProbPoly}), defined by (\ref{Ch4EqSolPoly+}) or (\ref{Ch4EqSolPoly*}), belongs to $\Li(\overline{\Omega^{+}}, \mu)$ and $\Li(\overline{\Omega^{-}},\mu)$ for $\mu < 1 - \frac{(n+1)(1-\nu)}{\mathfrak{m}(\mathcal{S})}$. 
	
	Now, let us suppose that there exist two solutions $\Phi_{1}$ and $\Phi_{2}$ to the problem (\ref{Ch4JumpProbPoly}), and define $\Phi := \Phi_{2} - \Phi_{1} $. This function is a solution to the homogeneous problem
	\begin{equation}\label{Ch4HomoJumpProbPoly}
		\begin{array}{ccc}
			(\mathcal{D}^{i}\Phi(x))^{+} - (\mathcal{D}^{i}\Phi(x))^{-} = 0 & x \in \mathcal{S} & 0 \leq i \leq k-1\\
			(\mathcal{D}^{i}\Phi(\infty))^{-} = 0 &  & 0 \leq i \leq k-1.
		\end{array}
	\end{equation}
	We shall prove that $\Phi \equiv 0$ is the unique solution to this problem such that $\mathcal{D}^{i}\Phi$  belongs to the classes $\Li(\overline{\Omega^{+}}, \mu)$ and $\Li(\overline{\Omega^{-}},\mu)$, for $i = 0,1,\cdots,k-1$. The proof is carried out by induction on $k$, by a repeated application of \cite[Theorem 10]{TamayoCastro2022}.
	
	Now we assume that (\ref{Ch4JumpProbPoly}) has the unique solution $\Phi \equiv 0$ such that $\mathcal{D}^{i}\Phi$  belongs to the classes $\Li(\overline{\Omega^{+}}, \mu)$ and $\Li(\overline{\Omega^{-}},\mu)$ for $i = 0,1,\cdots,k-1$; for $k = l$, and let us consider the problem for $k = l + 1$
	\begin{equation}\label{Ch4HomoJumpProbPolyT}
		\begin{array}{ccc}
			(\mathcal{D}^{i}\Phi(x))^{+} - (\mathcal{D}^{i}\Phi(x))^{-} = 0 & x \in \mathcal{S} & 0 \leq i \leq l \\
			(\mathcal{D}^{i}\Phi(\infty))^{-} = 0 &  & 0 \leq i \leq l.
		\end{array}
	\end{equation}
	Let $\Phi$ be a solution of (\ref{Ch4HomoJumpProbPolyT}). If we denote $\Psi := \mathcal{D}\Phi$, then $\mathcal{D}^{l}\Psi := \mathcal{D}^{l + 1}\Phi = 0$ in $\mathbb{R}^{n+1}\setminus\mathcal{S}$ and 
	\begin{equation*}
		\begin{array}{ccc}
			(\mathcal{D}^{i}\Psi(x))^{+} - (\mathcal{D}^{i}\Psi(x))^{-} = 0 & x \in \mathcal{S} & 0 \leq i \leq l - 1 \\
			(\mathcal{D}^{i}\Psi(\infty))^{-} = 0 &  & 0 \leq i \leq l - 1.
		\end{array}
	\end{equation*}
	Consequently, $\Psi$ represents a solution of (\ref{Ch4HomoJumpProbPoly}) with $k =l$. Then, by the induction hypothesis, $\Psi \equiv 0$ is the only solution in this class. As a result $\mathcal{D}\Phi = 0 $ in $\mathbb{R}^{n+1}\setminus\mathcal{S}$, and
	\begin{equation*}
		\begin{array}{ccc}
			\Phi(x)^{+} - \Phi(x)^{-} = 0 & x \in \mathcal{S}  \\
			\Phi(\infty)^{-} = 0 &.
		\end{array}
	\end{equation*}
	Therefore, as in the proof of \cite[Theorem 10]{TamayoCastro2022}, we have $\Phi \equiv 0$ in $\mathbb{R}^{n+1}$, and the proof is complete.
\end{proof}

\section{A Class of hypersurfaces in $\mathbb{R}^{n + 1}$}\label{Example}
A new solvability condition for RBVPs for polymonogenic functions via a high-dimensional Marcinkiewicz exponent has been proved. However, there naturally arises the question of whether this condition improves those involving the Minkowski dimension presented in \cite{AAB2015}? As a matter of fact, Theorem \ref{TheoIneqRm} is sufficient to guarantee that new condition does not ever be worth the effort of formulating them based on the Minkowski dimension. Indeed, a class of hypersurfaces in three-dimensional spaces so constructed in \cite{TamayoCastro2022} shows that strict inequality can occur in Theorem \ref{TheoIneqRm}, what is an extension to the case of the complex plane given in \cite{1DKats16}.

In this section, we generalize these constructions to a class of hypersurfaces in the $(n+1)$-dimensional space. Summarizing we have. 
\begin{theorem}\label{TheoInEx}
	Let $\alpha \geq 1$ and $\beta \geq n$. For each value $d \in (n, n + 1)$, there exists a uncountable collection of topologically compact hypersurfaces $\mathcal{S}^{n + 1}_{\alpha, \beta}$, which are the boundary of a Jordan domain in $\mathbb{R}^{n + 1}$ such that $d = \overline{\dim}_{M}(\mathcal{S}^{n + 1}_{\alpha, \beta})$ and $\mathfrak{m}(\mathcal{S}^{n + 1}_{\alpha, \beta}) > (n + 1) - \overline{\dim}_{M}(\mathcal{S}^{n + 1}_{\alpha, \beta})$ for suitable values of $\alpha$ and $\beta$.
\end{theorem}
\begin{proof}
	The proof consists in the construction of hypersurfaces $\mathcal{S}^{n+1}_{\alpha, \beta}$ having the desired properties.
	
	Let $ Q = [0, 1]\times[0, 1]\times[0, 1]\times\cdots\times[-1, 0]$ be a $(n+1)$-dimensional cube. We will add infinitely many $(n+1)$-dimensional rectangles with suitable dimensions to this cube. In order to do that, let us fix $\alpha \geq 1$ and $\beta \geq n$. Initially, we will conveniently divide the segment $[0,1]$ in the $x_{0}$ axis. We break down it into the sub-segments $[2^{-m},2^{-m+1}]$ for every $m \in \mathbb{N}$, and we divide each of this sub-segments into $2^{[m\beta]}$ equally spaced segments where $[m\beta]$ is the integer part of $m\beta$. The endpoints at the right side of these segments will be denoted by $y_{mj}$, where $j = 1, 2, ..., 2^{[m\beta]}$. In addition, let $a_{m} = 2^{-m - [m\beta]}$ be the distance between two consecutive points $y_{mj}$ and $y_{m(j+1)}$, and $C_{m} =\frac{1}{2}a_{m}^{\alpha}$. Then, let $R_{mj}$ be the $(n+1)$-dimensional rectangles defined as
	\begin{equation*}
		R_{mj} =  [y_{mj} - C_{m}, y_{mj}]\times[0, 2^{-m}]\times\cdots\times[0, 2^{-m}].
	\end{equation*}
	Hence, we define 
	\begin{equation*}
		T^{n+1}_{\alpha, \beta} := Q\bigcup\left(\bigcup_{m = 1}^{\infty}\bigcup_{j = 1}^{2^{[m\beta]}}R_{mj}\right).	
	\end{equation*}
	The claimed hypersurfaces $\mathcal{S}^{n+1}_{\alpha, \beta}$ are the boundaries of the corresponding $T^{n+1}_{\alpha, \beta}$. 
	
	We first obtain a suitable lower bound on the Marcinkiewicz exponent of hypersurfaces $\mathcal{S}^{n + 1}_{\alpha, \beta}$. To do this, split up the integral as follows
	\begin{equation}\label{IntegralMarcProf}
		\int\limits_{\Omega^{+}}\dfrac{dV}{\dist^{p}(x,\mathcal{S}^{n + 1}_{\alpha, \beta})} = \int\limits_{Q}\dfrac{dV}{\dist^{p}(x, \mathcal{S}^{n + 1}_{\alpha,\beta})} + \sum_{m = 1}^{\infty}\sum_{j = 1}^{2^{[m\beta]}} \int\limits_{R_{mj}}\dfrac{dV}{\dist^{p}(x, \mathcal{S}^{n + 1}_{\alpha,\beta})}.
	\end{equation}
	Let $L^{''}_{mj}$ be the subset of $R_{mj}$ such that $\dist(x, \mathcal{S}_{\alpha,\beta}) = \dist(x, \Gamma^{''}_{mj})$ where $\Gamma^{''}_{mj} = \left\lbrace y_{mj}\right\rbrace \times[0, 2^{-m}]\times\cdots\times[0, 2^{-m}]$, then
	\begin{equation*}
		\int\limits_{R_{mj}}\dfrac{dV}{\dist^{p}(x, \mathcal{S}^{n + 1}_{\alpha,\beta})} = \int\limits_{R_{mj}\setminus L^{''}_{mj}}\dfrac{dV}{\dist^{p}(x, \mathcal{S}^{n + 1}_{\alpha,\beta})} + \int\limits_{L^{''}_{mj}}\dfrac{dV}{\dist^{p}(x, \mathcal{S}^{n + 1}_{\alpha,\beta})}.
	\end{equation*}
	and 
	\begin{eqnarray*}
		\int\limits_{L^{''}_{mj}}\dfrac{dV}{\dist^{p}(x, \mathcal{S}^{n + 1}_{\alpha,\beta})} = \int\limits_{L^{''}_{mj}}\dfrac{dV}{\dist^{p}(x, \Gamma^{''}_{mj})}.
	\end{eqnarray*}
	
	Hence, using the estimate
	\begin{eqnarray*}
		\int\limits_{L^{''}_{mj}}\dfrac{dV}{\dist^{p}(x, \Gamma^{''}_{mj})} = \int\limits_{L^{''}_{mj}}\dfrac{dV}{\arrowvert y_{mj} - x_{0}\arrowvert^{p}} \leq \int\limits_{R_{mj}}\dfrac{dV}{\arrowvert y_{mj} - x_{1}\arrowvert^{p}}.
	\end{eqnarray*}
	and after some straightforward computations, we obtain that if the sum
	\begin{equation}\label{GeoSumProof}
		\sum_{m = 1}^{\infty}2^{m\beta - nm - (1 - p)\alpha(m + m\beta)},
	\end{equation}
	converges, then so is the integral in \eqref{IntegralMarcProf}. This geometric sum \eqref{GeoSumProof} converges if and only if the condition
	\begin{equation*}
		p < 1 - \dfrac{\beta - n}{\alpha(\beta+1)},
	\end{equation*}
	is fulfilled. Thus, the following estimation for the inner Marcinkiewicz exponent holds 
	\begin{equation*}
		\mathfrak{m}^{+}(\mathcal{S}^{n + 1}_{\alpha, \beta}):= \sup\{p>0: I_{p}(G^{+})<\infty\} \geq 1 - \dfrac{\beta - n}{\alpha(\beta+1)}.
	\end{equation*}
	We have proved more, namely that the absolute Marcinkiewicz exponent satisfies
	\begin{equation*}
		\mathfrak{m}(\mathcal{S}^{n + 1}_{\alpha, \beta}) := \max\{\mathfrak{m}^{+}(\mathcal{S}^{n + 1}_{\alpha, \beta}), \mathfrak{m}^{-}(\mathcal{S}^{n + 1}_{\alpha, \beta})\} \geq \mathfrak{m}^{+}(\mathcal{S}^{n + 1}_{\alpha, \beta}) \geq 1 - \dfrac{\beta - n}{\alpha(\beta+1)}.
	\end{equation*}
	We proceed to calculate the value of the Minkowski dimension of hypersurfaces $\mathcal{S}^{n + 1}_{\alpha, \beta}$. It is sufficient to show that its lower and upper bounds are the same.
	
	First, look for a suitable upper bound for $\overline{\dim}_{M}(\mathcal{S}^{n + 1}_{\alpha, \beta})$. We define the sets $\Lambda_{m} := \bigcup\limits_{j = 1}^{2^{[m\beta]}}[\partial R_{mj}\setminus(\partial R_{mj})\arrowvert_{x_{n} = 0}]$ and $\Lambda := \bigcup\limits_{m = 1}^{\infty}\Lambda_{m}$. Furthermore, $\widehat{Q} := \partial Q\setminus[\bigcup\limits_{m = 1}^{\infty}\bigcup\limits_{j = 1}^{2^{[m\beta]}}(\partial R_{mj})\arrowvert_{x_{n} = 0}]$, we can see that $\mathcal{S}^{n + 1}_{\alpha, \beta} = \widehat{Q}\cup\Lambda$.
	
	The following step makes use of the grid $\mathcal{M}_{k}$ defined in Subsection \ref{SubSecFracDim}. Initially, note that with $2(n+1)\left( \frac{1}{2^{-k}}\right) ^{n}$ cubes of $\mathcal{M}_{k}$, we can cover $\widehat{Q}$. In order to study $\Lambda$, we need to consider three cases. The first case is if $m < k$ and $C_{m} > 2^{-k}$, then $2^{[m\beta] +1}\left( \frac{2^{-m}}{2^{-k}}\right) ^{n}$ cubes have to be used to cover the faces of the $(n + 1)$-dimensional rectangles $R_{mj}$ parallel to $x_{0} = 0$, in $\Lambda_{m}$. No more than $2\left( \frac{2^{-m}}{2^{-k}}\right) ^{n}$ cubes are required to cover the $n$-dimensional rectangles in $\Lambda_{m}$ parallel to the coordinate plane $x_{l} = 0$, for every $l = 1, 2, \cdots, n$.
	
	The second, and main case, is if $m<k$ and $C_{m} \leq 2^{-k}$. Here we must study two more cases. When $C_{m} \leq 2^{-k}$, $k > m$, and also $a_{m} - C_{m} > 2^{-k}$, thus analogously than in the previous step $2^{[m\beta]}(\frac{2^{-m}}{2^{-k}})^{n}$ cubes are sufficient to cover the faces of $R_{mj}$'s parallel to $x_{0} = 0$ in $\Lambda_{m}$. While, to cover the $n$-dimensional rectangles in $\Lambda_{m}$ parallel to the coordinate plane $x_{l} = 0$ are not required more than $2\left( \frac{2^{-m}}{2^{-k}}\right) ^{n}$ cubes of $\mathcal{M}_{k}$, for every $l = 1, 2, \cdots, n$\\ 
	If $C_{m} \leq 2^{-k}$, $k > m$, and $a_{m} - C_{m} \leq 2^{-k}$, hence $\left( \frac{2^{-m}}{2^{-k}}\right) ^{n+1}$ cubes in $\mathcal{M}_{k}$ are adequate to cover $\Lambda_{m}$.
	
	Finally, when $m \geq k$, by definition, the hypersurfaces $\Lambda_{m}$, with $m > k$, are all covered by only one cube of the grid $\mathcal{M}_{k}$. While the hypersurface $\Lambda_{k}$ is covered by another of these cubes. Therefore, we obtain
	\begin{eqnarray*}
		N_{k}(\mathcal{S}^{n+1}_{\alpha,\beta}) \leq 2 + 2(n+1)\cdotp2^{nk} + 2\sum\limits_{C_{m}> 2^{-k},\,\,\, k>m}2^{[m\beta]+nk-nm} + 2n\sum\limits_{C_{m}> 2^{-k},\,\,\, k>m}2^{nk-nm} + \\
		+ \sum\limits_{C_{m}\leq 2^{-k},\,\,\, a_{m} - C_{m} \leq 2^{-k},\,\,\, k>m}2^{(n + 1)k-(n + 1)m} + \sum\limits_{C_{m}\leq 2^{-k} < a_{m} - C_{m},\,\,\, k>m}2^{[m\beta]+nk-nm} + \\
	\end{eqnarray*}
	\begin{equation*}
		+ 2n\sum\limits_{C_{m}\leq 2^{-k} < a_{m} - C_{m} \,\,\,	k>m}2^{nk-nm}. 
	\end{equation*}
	Working with the conditions on the sums in the previous inequality, we are able to obtain the following greater estimates
	\begin{equation*}
		N_{k}(\mathcal{S}^{n+1}_{\alpha,\beta}) \leq 2 + 2(n+1)\cdotp2^{nk} + 3\sum\limits_{2^{-k} < a_{m},\,\,\, k>m}2^{[m\beta]+nk-nm} + 4n\sum\limits_{2^{-k} < a_{m},\,\,\, k>m}2^{nk-nm} + 
	\end{equation*}
	\begin{equation*}
		+ \sum\limits_{\frac{a_{m}}{2} \leq 2^{-k},\,\,\, k>m}2^{(n+1)k-(n+1)m}. 
	\end{equation*}
	We will denote by $B_{k}$ and $H_{k}$ the integers defined by the conditions 
	\begin{eqnarray}\label{IneBk}
		\dfrac{k}{1+\beta}-1 \leq B_{k} < \dfrac{k}{1+\beta},\\
		\dfrac{k-1}{1+\beta}-1 \leq H_{k} < \dfrac{k-1}{1+\beta}.
	\end{eqnarray} 
	It is easy to check that the condition $a_{m} > 2^{-k}$ is satisfied if and only if $m = 1, 2,..., B_{k}$. By taking into account $H_{k}$ for the sum under the conditions $\frac{a_{m}}{2} \leq 2^{-k}, k>m$; and $B_{k}$ for those under the conditions $2^{-k} < a_{m}, k>m$ , we get through some estimates the following inequality
	\begin{equation*}
		N_{k}(\mathcal{S}^{n + 1}_{\alpha,\beta}) \leq D(k)2^{\frac{(n+1)k\beta}{\beta+1}},
	\end{equation*}
	where $D(k) = ak + c$; here $a$ and $c$ only depend on $\beta$ and $n$.
	Then,
	\begin{equation*}
		\overline{\dim}_{M}(\mathcal{S}^{n + 1}_{\alpha, \beta}) \leq \dfrac{(n+1)\beta}{\beta+1}.
	\end{equation*}
	Now we will compute a lower bound. In order to do that, we will build a set $A_{\beta}$ such that $A_{\beta} \subset \mathcal{S}^{n + 1}_{\alpha, \beta}$ and therefore $ \overline{\dim}_{M}(A_{\beta}) \leq \overline{\dim}_{M}(\mathcal{S}^{n + 1}_{\alpha, \beta})$.\\
	We will present the auxiliary $n$-dimensional rectangles $P_{mj}$
	\begin{equation*}
		P_{mj} = \{y_{mj}\}\times[0, 2^{-m}]\times\cdots\times[0, 2^{-m}],
	\end{equation*}
	and the set $A_{\beta}$ is given by
	\begin{equation*}
		A_{\beta} = \bigcup_{m = 1}^{\infty}\bigcup_{j = 1}^{2^{[m\beta]}}P_{mj}.
	\end{equation*}
	Observe, by construction, that $A_{\beta} \subset \mathcal{S}_{\alpha, \beta}$. The task is now to find a lower bound on $\overline{\dim}_{M}(A_{\beta})$. 
	
	The distance between $P_{mj}$ and $P_{mj+1}$ is $a_{m} = 2^{-m -[m\beta]}$. If $k > m$, and $a_{m} > 2^{-k}$, then two of these rectangles can not be intersected by the same cube in the grid $\mathcal{M}_{k}$. Then, $(\frac{2^{-m}}{2^{-k}})^{n}$  cubes in $\mathcal{M}_{k}$ cover a single $n$-dimensional rectangle $P_{mj}$, due to the fact that the lengths of each edge in these $n$-dimensional rectangles is $2^{-m}$.
	
	Therefore, the total number cubes required to cover the $2^{[m\beta]}$ rectangles $P_{mj}$ for a fixed $m$ is $2^{[m\beta]}(\frac{2^{-m}}{2^{-k}})^{n}$. Thus we get
	\begin{equation*}
		N_{k}(A_{\beta}) \geq 2\cdotp\sum\limits_{a_{m}> 2^{-k},\,\,\, k>m}2^{[m\beta]+nk-nm},
	\end{equation*}
	where $N_{k}(A_{\beta})$ is the smallest amount of cubes in $\mathcal{M}_{k}$ which cover $A_{\beta}$.\\
	Being $B_{k}$ as in \eqref{IneBk} we get
	\begin{equation*}
		\sum\limits_{a_{m}> 2^{-k},\,\,\, k>m}2^{[m\beta]+nk-nm} = 2^{nk}\sum\limits_{n = 1}^{B_{k}}2^{[m\beta]-nm} \geq 2^{nk-1}\sum\limits_{n = 1}^{B_{k}}2^{m(\beta-n)} \geq C2^{\frac{(n+1)k\beta}{\beta+1}},
	\end{equation*}
	where $C$ does not depend on $k$. Hence
	\begin{equation*}
		\overline{\dim}_{M}(S^{n+1}_{\alpha,\beta}) \geq \overline{\dim}_{M}(A_{\beta}) \geq \frac{(n+1)\beta}{\beta + 1}.
	\end{equation*}
	Therefore,
	\begin{equation*}
		\overline{\dim}_{M}(S^{n+1}_{\alpha,\beta}) = \frac{(n+1)\beta}{\beta + 1}.
	\end{equation*}
	Having disposed the value of the Minkowski dimension and a lower bound on the Marcinkiewicz exponent of every hypersurface $S^{n+1}_{\alpha,\beta}$, we are in a position to finish the proof of Theorem \ref{TheoInEx}. 
	
	Indeed, when $\alpha > 1 $ and $\beta > n$ we have
	\begin{equation*}
		\mathfrak{m}(\mathcal{S}_{\alpha, \beta}) \geq \mathfrak{m}^{+}(\mathcal{S}_{\alpha, \beta}) \geq 1 - \dfrac{\beta - n}{\alpha(\beta+1)} > 1 - \dfrac{\beta - n}{\beta+1} = (n + 1) - \dfrac{(n + 1)\beta}{\beta + 1} = (n + 1) - \overline{\dim}_{M}(\mathcal{S}_{\alpha, \beta}).
	\end{equation*}
	Setting $\beta = \displaystyle\frac{d}{(n + 1) - d}$, for every $d \in (n, n + 1)$, we obtain $\overline{\dim}_{M}(\mathcal{S}_{\alpha, \beta}) = d$ for each $\alpha > 1$. This means that $\{S^{n+1}_{\alpha,\beta}\}$ is an uncountable family, which completes the proof. 
\end{proof}

\textbf{Funding:} C. D. Tamayo Castro gratefully acknowledges the financial support of the Postgraduate Study Fellowship of the Consejo Nacional de Ciencia y Tecnología (CONACYT) (Grant Number 957110). J. Bory Reyes was partially supported by Instituto Politécnico Nacional in the framework of SIP programs
(SIP20230312).

\bibliographystyle{plain}
\bibliography{BiblioMarcPolyPreprint}

\end{document}